\documentclass[12pt,preprint]{elsarticle}

\usepackage{graphics,graphicx}
\usepackage{subfigure}
\usepackage{amssymb, amsfonts, amsmath,mathrsfs} 
\usepackage{eucal}

\usepackage{multirow}
\usepackage{lscape}
\usepackage{setspace}

\newcommand{\Eref}[1]{Equation (\ref{#1})}

\newcommand{\fref}[1]{Figure~\ref{#1}}

\newcommand{\BB}{\mathbf{B}}

\newcommand{\DD}{\mathbf{D_b}}

\newcommand{\KK}{\mathbf{K}}

\newcommand{\bn}{\mathbf{N}}

\newcommand{\xx}{\mathbf{x}}

\newcommand{\bfm}{\mathbf{M}}

\newcommand{\rmd}{\mathrm{d}}

\newcommand{\bveps}{\boldsymbol{\varepsilon}}
\newcommand{\bvsig}{\boldsymbol{\sigma}}

\journal{Computational Material Science}

\begin{document}

\begin{frontmatter}
\title{Free flexural vibration of functionally graded size-dependent nanoplates}

\author[a]{S Natarajan\fnref{label2}\corref{cor1}}
\author[b]{S Chakraborty}
\author[b]{M Thangavel}

\address[a]{Postdoctoral Research Fellow, Department of Aerospace Engineering, Indian Institute of Science, Bangalore, INDIA.}
\address[b]{Research Assistant, Department of Aerospace Engineering, Indian Institute of Science, Bangalore, INDIA.}

\fntext[3]{sundararajan.natarajan@gmail.com}

\begin{abstract}
In this paper, the linear free flexural vibration behaviour of functionally graded (FG) size-dependent nanoplates are investigated using the finite element method. The field variables are approximated by non-uniform rational B-splines. The size-dependent FG nanoplate is investigated by using Eringen's differential form of nonlocal elasticity theory. The material properties are assumed to vary only in the thickness direction and the effective properties for FG nanoplate are computed using Mori-Tanaka homogenization scheme. The accuracy of the present formulation is tested considering the problems for which solutions are available. A detailed numerical study is carried out to examine the effect of material gradient index, the characteristic internal length, the plate thickness, the plate aspect ratio and the boundary conditions on the global response of FG nanoplate.
\end{abstract}

\begin{keyword} 
Functionally graded, Mori-Tanaka, Eringen's gradient elasticity, partition of unity, finite element, NURBS, internal length.
\end{keyword}

\end{frontmatter}

\section{Introduction}

Rapid advancement in the application of micro/nano structures in engineering fields, mainly MEMS/NEMS devices, due to their superior mechanical properties, rendered a sudden momentum in modelling the structures of nano and micro length scale. It has been observed that there is a considerable difference in the structural behaviour of material at micro-/nano-scale when compared to their bulk counterpart. The difficulty in using the classical theory is that these classical models fails to capture the size effects. The classical model over predicts the response of micro-/nano-structures. Also, in the classical model, the particles influence one another by contact forces. Another way to capture the size effect is to rely on first principle calculations or molecular dynamic simulations (MD). But even the MD simulation at micron scale is computationally demanding. Hence, there was a need for some theory that can incorporate the length scale effect. Several modifications of the classical elasticity formulation have been proposed in the literature~\cite{mindlin1964,kroner1967,eringen1972,aifantis1984,askesmetrikine2008}. A common feature of all these theories is that they include one or several intrinsic length scales and the particles influence one another by long range cohesive forces. In this way, the internal length scale can be considered in the constitutive equations simply as a material parameter, also referred to as `nonlocal parameter' and `internal characteristic length' in the literature. \\

Among the various nonlocal theories, Eringen's differential gradient elasticity~\cite{eringen1972} has received considerable attention to study the static and the dynamic characteristics of nano beams~\cite{aydogdu2009a,hsulee2011} and plates~\cite{luzhang2007,murmupradhan2009a,reddypang2008}. The theory has been applied by several authors to study the axial vibrations~\cite{murmupradhan2009a,aydogdu2009,daneshfarajpour2012,filizaydogdu2010} and free transverse vibrations of nanostructures~\cite{ansariarash2011,aydogdu2009,loyal'opez-puente2009,murmupradhan2009}. Reddy and Pang~\cite{reddypang2008} derived governing equations of motion for Euler-Bernoulli beams and Timoshenko beams using the nonlocal differential relations of Eringen~\cite{eringen1972}. The classical plate theory (CPT) and first order shear deformation theory (FSDT) have been reformulated using the nonlocal differential constitutive relations in ~\cite{luzhang2007,reddypang2008,reddy2010}. It is observed that the increasing effect of the nonlocal parameter is to decrease the fundamental frequencies of the structure. In all the earlier studies, the nonlocal parameter is kept as a variable and effect of the nonlocal parameter on the structural response is studied. One approach is to estimate the nonlocal parameter by matching the phonon dispersion relations computed by these theories with the lattice dynamics dispersion relation~\cite{chenlee2003}. Ansari \textit{et al.,}~\cite{ansarisahmani2010} used MD simulations to derive the appropriate value of the nonlocal parameter. Recently, Eltaher~\textit{et al.,}~\cite{eltaheremam2012} employed Eringen's differential theory to study the size-effect on the structural response of functionally graded nanobeams. Chang~\cite{chang2012} studied the axial vibration characteristics of non-uniform and non-homogeneous nanorods. It was observed that the nonlocal parameter has profound impact on the structural response. \\

The existing numerical approaches are limited to using radial basis functions~\cite{roqueferreira2011} and differential quadrature method~\cite{murmupradhan2009a,eltaheremam2012, chang2012} to study the response of nanostructures based on nonlocal elasticity theory. The main objective of this paper is to employ Eringen's differential gradient elasticity to study the structural response of functionally graded (FG) size dependent nanobeams. The plate kinematics is described by FSDT and the governing equations of motion are solved within a finite element (FE) framework. The field variables are approximated by non-uniform rational B-splines (NURBS)~\cite{piegltiller1996}. The effect of the nonlocal parameter, the material gradient index, the plate aspect ratio, thickness of the plate and boundary conditions on the fundamental frequencies are numerically studied.\\

The paper is orgranized as follows. A brief over of FG material, Eringen's nonlocal elasticity and Reissner-Mindlin plate theory is given in the next section. Section \ref{spatdis} presents an overview of NURBS basis functions and spatial discretization. The numerical results and parametric studies are given in Section \ref{numresults}, followed by concluding remarks in the last section.

%------------ Theoretical formulation
\section{Theoretical development}
\label{fgmsyn}

%------------------- Functionally graded material
%----------- Functionally graded material
\subsection{Functionally graded material}
\begin{figure}[htpb]
\centering
\input{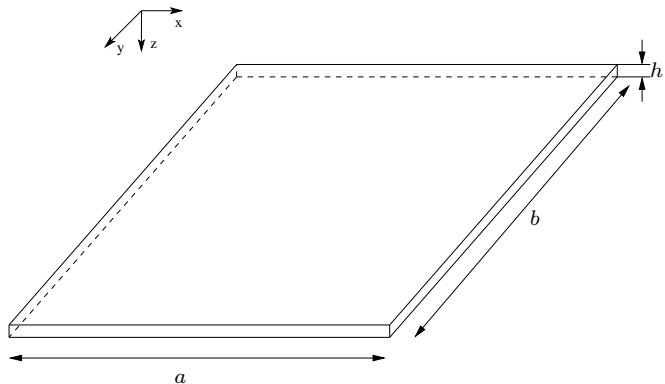}
\caption{Coordinate system of a rectangular FGM plate.}
\label{fig:platefig}
\end{figure}

A rectangular plate (length $a$, width $b$ and thickness $h$) with functionally graded material (FGM), made by mixing two distinct material phases: a metal and ceramic is considered with coordinates $x,y$ along the in-plane directions and $z$ along the thickness direction (see \fref{fig:platefig}). The material on the top surface $(z=h/2)$ of the plate is ceramic and is graded to metal at the bottom surface of the plate $(z=-h/2)$ by a power law distribution. The homogenized material properties are computed using the Mori-Tanaka Scheme~\cite{moritanaka1973,benvensite1987}. 

Based on the Mori-Tanaka homogenization method, the effective bulk modulus $K$ and shear modulus $G$ of the FGM are evaluated as~\cite{moritanaka1973,benvensite1987,chengbatra2000,qianbatra2004}:

\begin{eqnarray}
{K - K_m \over K_c - K_m} &=& {V_c \over 1+(1-V_c){3(K_c - K_m) \over 3K_m + 4G_m}} \nonumber \\
{G - G_m \over G_c - G_m} &=& {V_c \over 1+(1-V_c){(G_c - G_m) \over G_m + f_1}}
\label{eqn:bulkshearmodulus}
\end{eqnarray}

where,

\begin{equation}
f_1 = {G_m (9K_m + 8G_m) \over 6(K_m + 2G_m)}
\end{equation}

Here, $V_i~(i=c,m)$ is the volume fraction of the phase material. The subscripts $c$ and $m$ refer to the ceramic and metal phases, respectively. The volume fractions of the ceramic and metal phases are related by $V_c + V_m = 1$, and $V_c$ is expressed as:

\begin{equation}
V_c(z) = \left( {2z + h \over 2h} \right)^n, \hspace{0.2cm}  n \ge 0
\label{eqn:volFrac}
\end{equation}

where $n$ in \Eref{eqn:volFrac} is the volume fraction exponent, also referred to as the gradient index. \fref{fig:volfrac} shows the variation of the volume fractions of ceramic and metal, respectively, in the thickness direction $z$ for the FGM plate. The effective Young's modulus $E$ and Poisson's ratio $\nu$ can be computed from the following expressions:

\begin{figure}[htpb]
\includegraphics[scale=1]{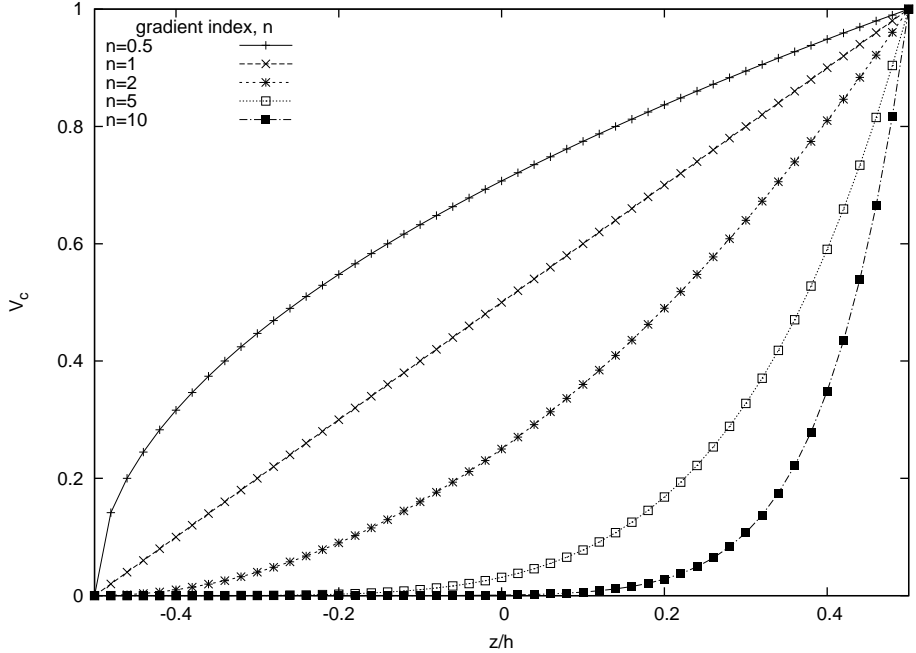}
\caption{Through thickness variation of volume fraction}
\label{fig:volfrac}
\end{figure}

\begin{eqnarray}
E = {9KG \over 3K+G} \nonumber \\
\nu = {3K - 2G \over 2(3K+G)}
\label{eqn:young}
\end{eqnarray}

The effective mass density $\rho$ is given by the rule of mixtures as~\cite{natarajanbaiz2011}:

\begin{equation}
\rho = \rho_c V_c + \rho_m V_m
\label{eqn:mdensity}
\end{equation}

%----------- Nonlocal elasticity
\subsection{Overview of nonlocal elasticity}
Eringen~\cite{eringen1972,eringenedelen1972}, by including the long range cohesive forces, proposed a nonlocal theory in which the stress at a point depends on the strains of an extended region around that point in the body. Thus, the nonlocal stress tensor $\bvsig$ at a point $\xx$ is expressed as:

\begin{equation}
\bvsig = \int_\Omega \alpha(|\xx^\prime - \xx|,\tau) \mathbf{t}(\xx^\prime)~\rmd \xx^\prime 
\label{eqn:eringenInte}
\end{equation}

where $\mathbf{t}(\xx)$ is the classical, macroscopic stress tensor at a point $\xx$ and the kernel function $\alpha(|\xx^\prime - \xx|)$ is called the nonlocal modulus, also referred to as the attenuation or influence function. The nonlocality effects at a reference point $\xx$ produced by the strains at $\xx$ and $\xx^\prime$ are included in the constitutive law by this function. The kernel weights the effects of the surrounding stress states. From the structure of the nonlocal constitutive equation, given by \Eref{eqn:eringenInte}, it can be seen that the nonlocal modulus has the dimension of (length)$^{-3}$~\cite{eringen1972,eringenedelen1972}. Typically, the kernel is a function of the Euclidean distance between the points $\xx$ and $\xx^\prime$. The nonlocal modulus has a very important role to play capturing the nonlocal effects. Ofcourse, ultimately, the accurate estimation of $\alpha(|\xx^\prime - \xx|)$ would dictate the overall reliability of a nonlocal model. Eringen~\cite{eringen1972,eringen2002} numerically determined the functional form of the kernel. By appropriate choice of the kernel function, Eringen~\cite{eringen1972} showed that the nonlocal constitutive equation given in integral form (see \Eref{eqn:eringenInte}) can be represented in an equivalent differential form as:

\begin{equation}
(1 - \tau^2L^2 \nabla^2) \bvsig = \mathbf{t}
\end{equation}

where $\tau^2 = \frac{\mu}{L^2} = \left( \frac{e_oa}{L} \right)^2$, $e_o$ is a material constant and $a$ and $L$ are the internal and external characteristic lengths, respectively.

%------- Reissner-Mindline plate
\subsection{Reissner-Mindlin plate}
Consider a rectangular Cartesian coordinate system $(x,y,z)$ with $xy-$ plane coinciding with the undeformed middle plane of the plate and the $z$- coordinate is taken positive downwards. The displacement field based on first order shear deformation plate theory (FSDT), also referred to as the Mindlin plate theory is given by:

\begin{eqnarray}
u(x,y,z,t) &=& u_o(x,y,t) + z \theta_x(x,y,t) \nonumber \\
v(x,y,z,t) &=& v_o(x,y,t) + z \theta_y(x,y,t) \nonumber \\
w(x,y,z,t) &=& w_o(x,y,t) 
\label{eqn:displacements}
\end{eqnarray}

\begin{equation}
\bveps  = \left\{ \begin{array}{c} \bveps_p \\ 0 \end{array} \right \}  + \left\{ \begin{array}{c} z \bveps_b \\ \bveps_s \end{array} \right\} 
\label{eqn:strain1}
\end{equation}

The midplane strains $\bveps_p$, bending strain $\bveps_b$, shear strain $\bveps_s$ in \Eref{eqn:strain1} are written as

\begin{equation}
\renewcommand{\arraystretch}{1.5}
\bveps_p = \left\{ \begin{array}{c} u_{o,x} \\ v_{o,y} \\ u_{o,y}+v_{o,x} \end{array} \right\}, \hspace{0.5cm}
\renewcommand{\arraystretch}{1.5}
\bveps_b = \left\{ \begin{array}{c} \theta_{x,x} \\ \theta_{y,y} \\ \theta_{x,y}+\theta_{y,x} \end{array} \right\} \hspace{0.5cm}
\renewcommand{\arraystretch}{1.5}
\bveps_s = \left\{ \begin{array}{c} \theta _x + w_{o,x} \\ \theta _y + w_{o,y} \end{array} \right\}.
\renewcommand{\arraystretch}{1.5}
\end{equation}

where the subscript `comma' represents the partial derivative with respect to the spatial coordinate succeeding it. The equations of motion are given by:

\begin{eqnarray}
\frac{\partial N_{xx}}{\partial x} + \frac{\partial N_{xy}}{\partial x} = I_{11} \ddot{u} + I_{12} \ddot{\theta_x} \nonumber \\
\frac{\partial N_{xy}}{\partial x} + \frac{\partial N_{yy}}{\partial x} = I_{11} \ddot{v} + I_{12} \ddot{\theta_y} \nonumber \\
\frac{\partial N_{xz}}{\partial x} + \frac{\partial N_{yz}}{\partial x} = I_{11} \ddot{w} \nonumber \\
\frac{\partial M_{xx}}{\partial x} + \frac{\partial M_{xy}}{\partial x} = I_{12} \ddot{u} + I_{22} \ddot{\theta_x} \nonumber \\
\frac{\partial M_{xy}}{\partial x} + \frac{\partial M_{yy}}{\partial x} = I_{12} \ddot{v} + I_{22} \ddot{\theta_y}
\label{eqn:plateeqeqn}
\end{eqnarray}

where $(I_{11}, I_{12}, I_{22}) = \int\limits_{-h/2}^{h/2} \rho(1,z,z^2)~\rmd z$. The nonlocal membrane stress resultants $\bn$ and the bending stress resultants $\bfm$ can be related to the membrane strains, $\bveps_p$ and bending strains $\bveps_b$ through the following constitutive relations~\citep{reddy2010}:

\begin{eqnarray}
\bn &=& \left\{ \begin{array}{c} N_{xx} \\ N_{yy} \\ N_{xy} \end{array} \right\}(1- \mu \nabla^2) = \mathbf{A}_e \bveps_p + \BB_{be} \bveps_b \nonumber \\
\bfm &=& \left\{ \begin{array}{c} M_{xx} \\ M_{yy} \\ M_{xy} \end{array} \right\}(1- \mu \nabla^2) = \BB_{be} \bveps_p + \DD_b \bveps_b  
\label{eqn:platenonlocal}
\end{eqnarray}

where the matrices $\mathbf{A}_e = A_{ij}, \BB_{be} = B_{ij}$ and $\DD_b = D_{ij}; (i,j=1,2,6)$ are the extensional, bending-extensional coupling and bending stiffness coefficients and are defined as

\begin{equation}
\left\{ A_{ij}, ~B_{ij}, ~ D_{ij} \right\} = \int_{-h/2}^{h/2} \overline{Q}_{ij} \left\{1,~z,~z^2 \right\}~dz, \hspace{0.2cm} i=1,2,6
\end{equation}

Similarly, the transverse shear force $Q = \{Q_{xz},Q_{yz}\}$ is related to the transverse shear strains $\bveps_s$ through the following equation

\begin{equation}
Q_{ij} = E_{ij} \bveps_s
\end{equation}

where $E_{ij} = \mathbf{E} = \int_{-h/2}^{h/2} \overline{Q} \upsilon_i \upsilon_j~dz;~ (i,j=4,5)$ is the transverse shear stiffness coefficient, $\upsilon_i, \upsilon_j$ is the transverse shear correction factors for non-uniform shear strain distribution through the plate thickness. The stiffness coefficients $\overline{Q}_{ij}$ are defined as

\begin{eqnarray}
\overline{Q}_{11} = \overline{Q}_{22} = {E(z) \over 1-\nu^2}; \hspace{1cm} \overline{Q}_{12} = {\nu E(z) \over 1-\nu^2}; \hspace{1cm} \overline{Q}_{16} = \overline{Q}_{26} = 0 \nonumber \\
\overline{Q}_{44} = \overline{Q}_{55} = \overline{Q}_{66} = {E(z) \over 2(1+\nu) }
\end{eqnarray}

The equations of motion in terms of midplane displacements and rotations are obtained by substituting the nonlocal constitutive relations (\Eref{eqn:platenonlocal}) into the equations of motion~\Eref{eqn:plateeqeqn}. For detailed derivation, interested readers are referred to the literature~\citep{murmupradhan2009a,reddy2010,malekzadehsetoodeh2011}.

%------------ spatial discretization
\section{Spatial discretization}
\label{spatdis}
In this study, the finite element model has been developed using NURBS basis function. 

\paragraph{Non-uniform Rational B-Splines}
The key ingredients in the construction of NURBS basis functions are: the knot vector (a non decreasing sequence of parameter value, $\xi_i \le \xi_{i+1}, i = 0,1,\cdots,m-1$), the control points, $P_i$, the degree of the curve $p$ and the weight associated to a control point, $w$. The i$^{th}$ B-spline basis function of degree $p$, denoted by $N_{i,p}$ is defined as~\cite{piegltiller1996,hughescottrell2005}:

\begin{eqnarray}
N_{i,0}(\xi) = \left\{ \begin{array}{cc} 1 & \textup{if} \hspace{0.2cm} \xi_i \le \xi \le \xi_{i+1} \\
0 & \textup{else} \end{array} \right. \nonumber \\
N_{i,p}(\xi) = \frac{ \xi- \xi_i}{\xi_{i+p} - \xi_i} N_{i,p-1}(\xi) + \frac{\xi_{i+p+1} - \xi}{\xi_{i+p+1}-\xi_{i+1}}N_{i+1,p-1}(\xi)
\end{eqnarray}

A $p^{th}$ degree NURBS curve is defined as follows:

\begin{equation}
\mathbf{C}(\xi) = \frac{\sum\limits_{i=0}^m N_{i,p}(\xi)w_i \mathbf{P}_i}{\sum\limits_{i=0}^m N_{i,p}(\xi)w_i}
\end{equation}

\begin{figure}[htpb]
\centering
\includegraphics[scale=0.8]{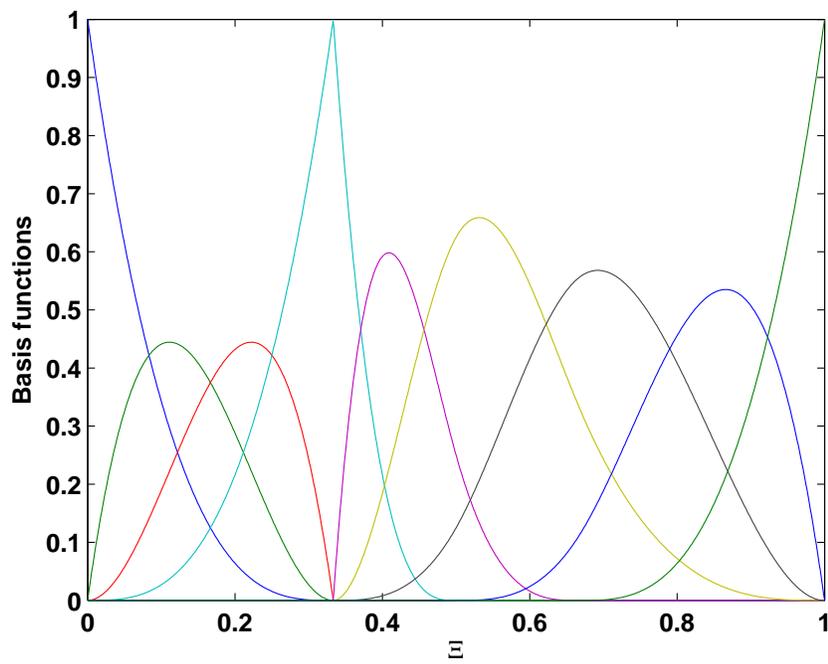}
\caption{non-uniform rational B-splines, order of the curve = 3}
\label{fig:nurbsplot}
\end{figure}

where $\mathbf{P}_i$ are the control points and $w_i$ are the associated weights. \fref{fig:nurbsplot} shows the third order non-uniform rational B-splines for a knot vector, $\Xi = \{0,~ 0,~ 0,~ 0,~ 1/3,~ 1/3,~ 1/3,~ 1/2,~ 2/3,~ 1,~ 1,~ 1,~ 1\}$. The NURBS basis functions has the following properties: (i) non-negativity, (ii) partition of unity, $\sum\limits_i N_{i,p} = 1$; (iii) interpolatory at the end points. These properties, make them best suitable to represent the displacement field in elasticity. As the same function is also used to represent the geometry, the exact representation of the geometry is preserved. It should be noted that the continuity of the NURBS functions can be custom tailored to suit the needs of the problem. Klassen~\textit{et al.,}~\cite{klassenmuller2010} and Fischer~\textit{et al.,}~\cite{fischerklassen2011} employed NURBS basis functions as approximating functions for problems in 2D gradient elasticity. The displacement field variable is approximated by the NURBS basis functions. By multiplying the governing equations of motion by appropriate test functions and by employing Bubnov-Galerkin method, we get the following set of algebraic equations:

\begin{equation}
\bfm \ddot{\boldsymbol{\delta}} + \KK \boldsymbol{\delta} = 0
\end{equation}

where $\bfm$ is the consistent mass matrix, $\KK$ is the stiffness matrix and $\boldsymbol{\delta}$ is the vector of the degree of freedom associated to the displacement field. After substituting the characteristic of the time function $\ddot{\boldsymbol{\delta}} = - \omega^2 \boldsymbol{\delta}$, the following algebraic equation is obtained:

\begin{equation}
( \KK - \omega^2 \bfm) \boldsymbol{\delta}  = \mathbf{0}
\end{equation}

where $\omega$ is the natural frequency. The frequencies are obtained using the standard generalized eigenvalue algorithm.

%------------- Numerical example section
\section{Numerical Example}
\label{numresults}
\begin{figure}[htpb]
\centering
%\subfigure[Simply supported]{\input{SimplyS.pstex_t}}
%\subfigure[clamped]{\input{Clamped.pstex_t}}
\input{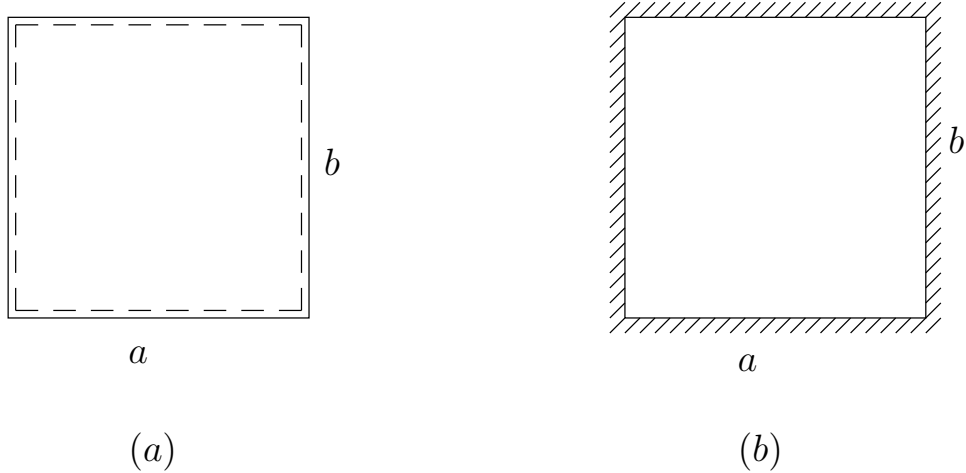}
\caption{Plate with boundary conditions: (a) all edges simply supported (SSSS) and (b) all edges clamped support (CCCC).}
\label{fig:simplyss}
\end{figure}

In this section, the free vibration characteristics of FG size-dependent nanoplate is numerically studied. \fref{fig:simplyss} shows the geometry and boundary conditions of the plate. The effect of gradient index $n$, internal characteristic length $\mu$, the plate aspect ratio $a/b$, thickness of the plate $h$ and boundary conditions on the natural frequency are numerically studied. The assumed values for the parameters are listed in Table \ref{table:parametervalues}. 

\begin{table}
\renewcommand\arraystretch{1.5}
\caption{Details of parameters used in the numerical study}
\centering
\begin{tabular}{ll}
\hline
Parameter & Assumed values  \\
\hline
Gradient index, $n$ & 0, 1, 2, 3, 4, 5 \\
Aspect ratio $a/b$ & 1, 2 \\
Thickness of the plate $a/h$ & 10, 20 \\
Internal length $\mu$ (nm$^2$) & 0, 1, 2, 3, 4, 5 \\
Boundary conditions & Simply supported and Clamped \\
\hline
\end{tabular}
\label{table:parametervalues}
\end{table}

In all cases, we present the non-dimensionalized free flexural frequency defined as:

\begin{equation}
\Omega = \omega h \sqrt{ \frac{ \rho_c}{G_c}}
\label{eqn:vibnondim}
\end{equation}

where $\omega$ is the natural frequency, $\rho_c$ and $G_c$ are the mass density and shear modulus of the ceramic, respectively. In the present numerical study, the effect of the nonlocal parameter is studied by the frequency ratio, defined as:

\begin{equation}
\textup{Frequency ratio} = \frac{\Omega_{NL}}{\Omega_L}
\end{equation}

where $\Omega_{NL}$ and $\Omega_L$ are the non-dimensionalized frequency based on nonlocal and local elasticity, respectively. The frequency ratio is a measure of the error made by neglecting the nonlocal effects. For example, if $\frac{\Omega_{NL}}{\Omega_L}=$ 1, the effect of the nonlocal parameter is not significant and becomes important for any other value.

The top surface is ceramic rich and the bottom surface is metal rich. The FG plate considered here is made up of silicon nitride (Si$_3$N$_4$) and stainless steel (SUS304). The mass density $\rho$ and the Young's modulus $E$ are: $\rho_c=$ 2370 kg/m$^3$, $E_c=$ 348.43e$^9$ N/m$^2$ for Si$_3$N$_4$ and $\rho_m=$ 8166 kg/m$^3$, $E_m=$ 201.04e$^9$ N/m$^2$ for SUS304. Poisson's ratio $\nu$ is assumed to be constant and taken as 0.3 for the current study. Here, the modified shear correction factor obtained based on energy equivalence principle as outlined in~\cite{singhprakash2011} is used.

The boundary conditions for simply supported and clamped cases are (see \fref{fig:simplyss}):

\noindent \emph{Simply supported boundary condition}: \\
\begin{eqnarray}
u_o = w_o = \theta_y = 0 \hspace{1cm} ~\textup{on} ~ x=0,a \nonumber \\
v_o = w_o = \theta_x = 0 \hspace{1cm} ~\textup{on} ~ y=0,b
\end{eqnarray}

\noindent \emph{Clamped boundary condition}: \\
\begin{eqnarray}
u_o = w_o = \theta_y = v_o = \theta_x = 0 \hspace{1cm} ~\textup{on} ~ x=0,a \nonumber \\
u_o = w_o = \theta_y = v_o = \theta_x = 0 \hspace{1cm} ~\textup{on} ~ y=0,b
\end{eqnarray}

\begin{figure}[htpb]
\centering
\includegraphics[scale=0.8]{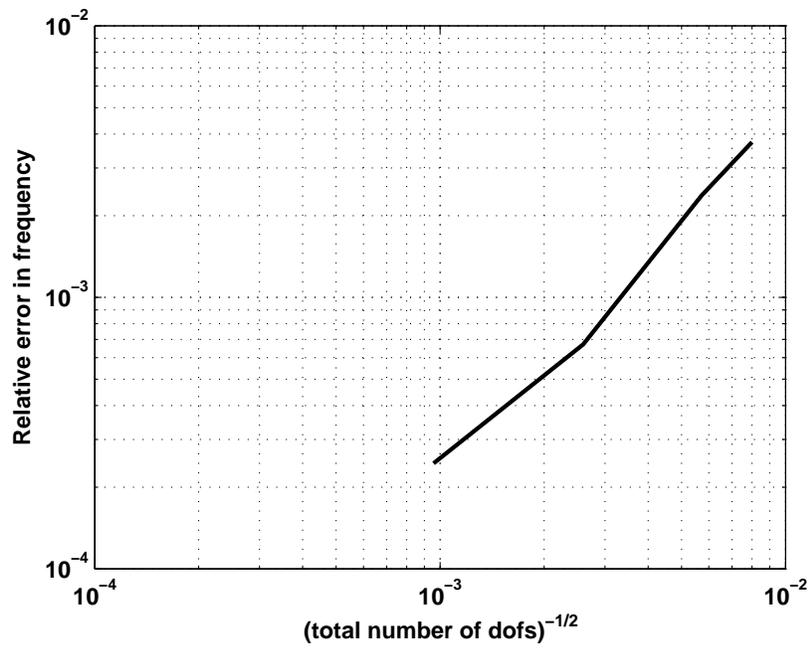}
\caption{The convergence of first fundamental frequency with total degrees of freedom. A 3$^{rd}$ order NURBS curve is used.}
\label{fig:meshcon}
\end{figure}

Before proceeding with the detailed study on the effect of gradient index and the characteristic internal length on the natural frequencies, the formulation developed is validated against the available results pertaining to the linear frequencies for an isotropic nanoplate based on FSDT~\cite{aghababaeireddy2009}. Table~\ref{table:SSisoplate_comp} gives the results of non-dimensional frequency for a simply supported square plate for different values of nonlocal parameter, the plate thickness and the plate aspect ratio. The numerical results from the present formulation are found to be in very good agreement with the existing solutions. It can be seen that the nonlocal theory predicts smaller values of natural frequency than the local elasticity theory.\fref{fig:meshcon} shows the convergence of the first fundamental frequency with increasing degrees of freedom (dofs). It can be seen that with increasing dofs, the analytical solution is approached and the method shows a constant convergence. 

% nondimensionalized frequencies
\begin{table}[htpb]
\renewcommand\arraystretch{1.5}
\caption{Comparison of fundamental frequency $\omega h \sqrt{\rho/G}$ for simply supported plate. $^{\dagger\dagger}$\cite{aghababaeireddy2009}, $^\sharp$Order of curve = 3, 5 $\times$ 5 control points.} 
\centering
\begin{tabular}{cccll}
\hline
$a/b$ & $a/h$ & $\mu$ & Ref.$^{\dagger\dagger}$ & NURBS$^\sharp$ \\
\hline
\multirow{6}{*}{1} & \multirow{3}{*}{10} & 0 & 0.0930 &  0.0929\\
& & 1 & 0.0850 &  0.0849\\
& & 5 & 0.0660 &  0.0659\\
\cline{3-5}
& \multirow{3}{*}{20} & 0 & 0.0239 & 0.0239\\
& & 1 & 0.0218 &  0.0219\\
& & 5 & 0.0169 &  0.0170\\
\hline
\multirow{6}{*}{2} & \multirow{3}{*}{10} & 0 & 0.0589 &  0.0590\\
& & 1 & 0.0556 &  0.0556\\
& & 5 & 0.0463 &  0.0464\\
\cline{3-5}
& \multirow{3}{*}{20} & 0 & 0.0150 & 0.0151\\
& & 1 & 0.0141 & 0.0141\\
& & 5 & 0.0118 & 0.0118\\
\hline
\end{tabular}
\label{table:SSisoplate_comp}
\end{table}

\begin{figure}[htpb]
\centering
\includegraphics{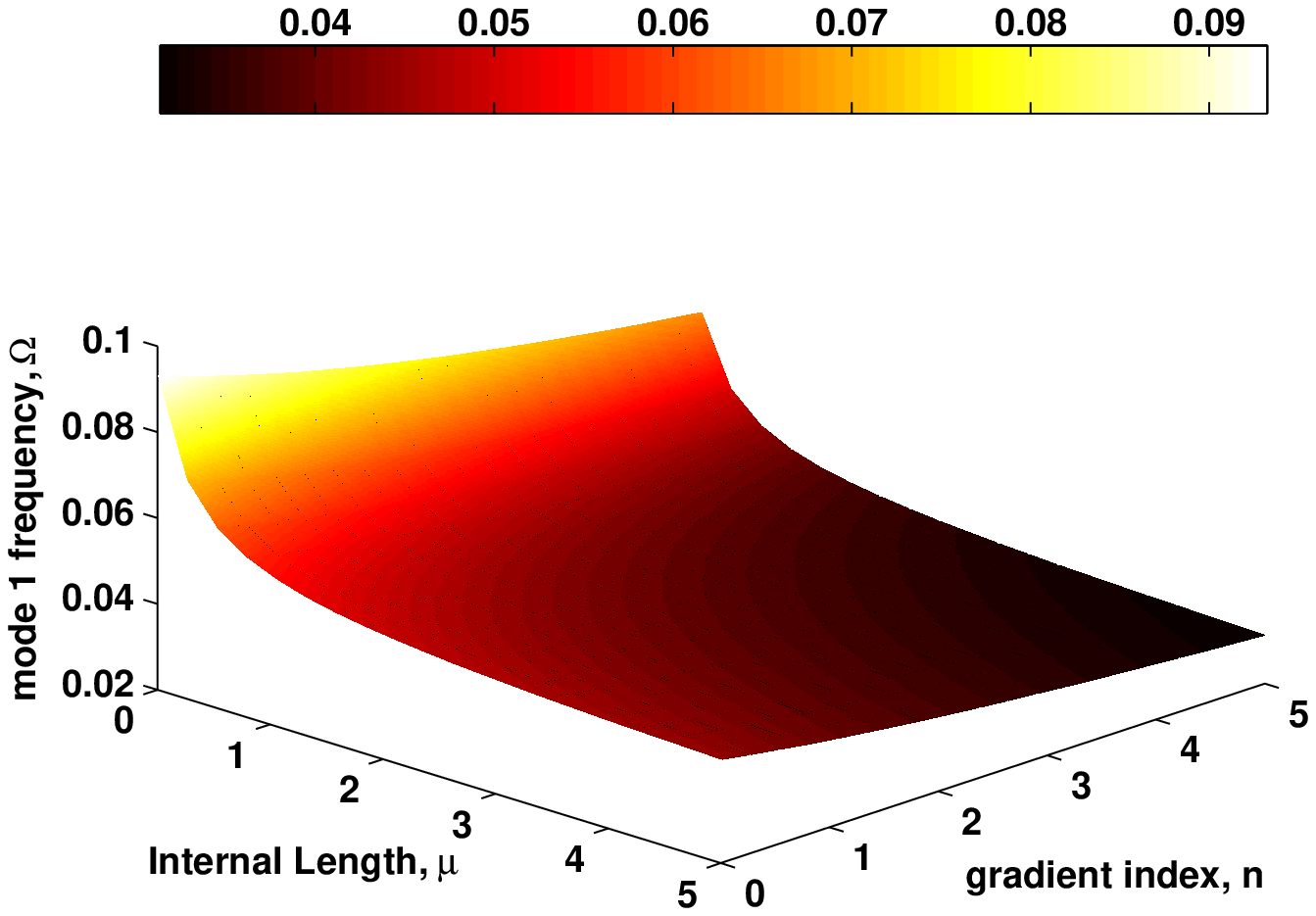}
\caption{Effect of characteristic internal length and gradient index on the first frequency for a simply supported square plate with $a/h=$ 10, $a/b=$ 1 and $a=$10 nm.}
\label{fig:mode1gIndIntLen}
\end{figure}

%------- numerical results
Next, the linear free flexural vibration of FG size-dependent nanoplate is numerically studied. The variation of the first natural frequency with the gradient index and the nonlocal parameter is graphically shown in \fref{fig:mode1gIndIntLen} for a simply supported square plate.

The values of the nonlocal parameter are assumed to vary between $\mu=$ 0 (local elasticity) to $\mu=$ 5 nm$^2$ and the gradient index is assumed to vary between $n=$ 0 (ceramic) to $n=$ 5. The combined effect of the nonlocal parameter $\mu$ and the gradient index $n$ is to reduce the natural frequency. The decrease in frequency due to increase in the gradient index is due to increase in the metallic volume fraction.

\begin{figure}[htpb]
\centering
\includegraphics[scale=0.8]{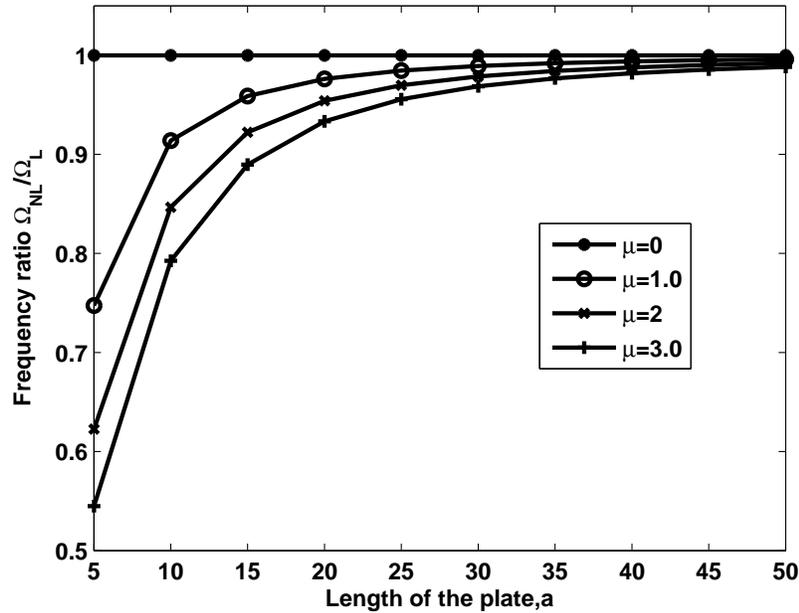}
\caption{Effect of plate dimension on the frequency ratio $(\Omega_{NL}/\Omega_L$ for a simply supported square plate for various internal length with $a/h=$ 10 and gradient index $n$= 5.}
\label{fig:platedimInfluence}
\end{figure}

\begin{figure}[htpb]
\centering
\includegraphics[scale=0.8]{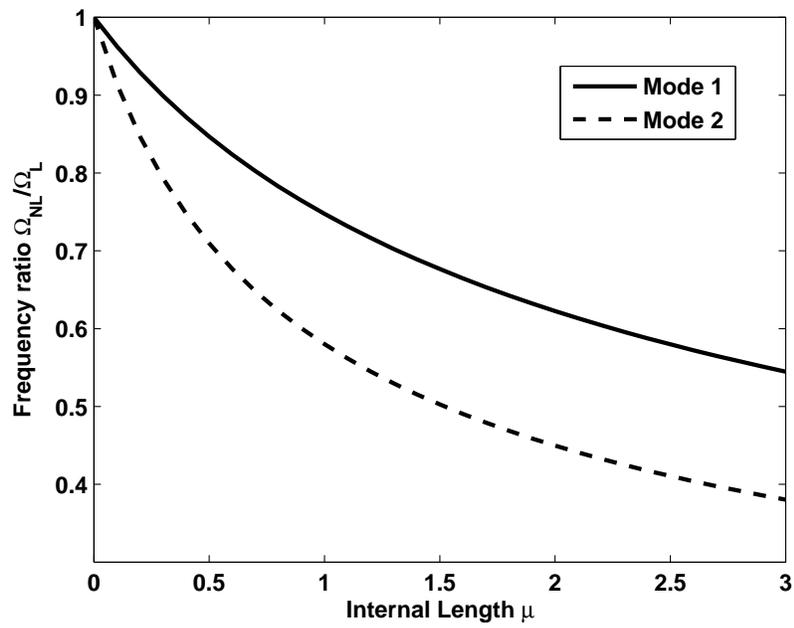}
\caption{Effect of the nonlocal parameter $\mu$ on the frequency ratio $(\Omega_{NL}/\Omega_L$ for a simply supported square plate for the first two frequencies with $a/h=$ 10 and gradient index $n=$ 5.}
\label{fig:diffmodes}
\end{figure}

\begin{figure}[htpb]
\centering
\includegraphics[scale=0.8]{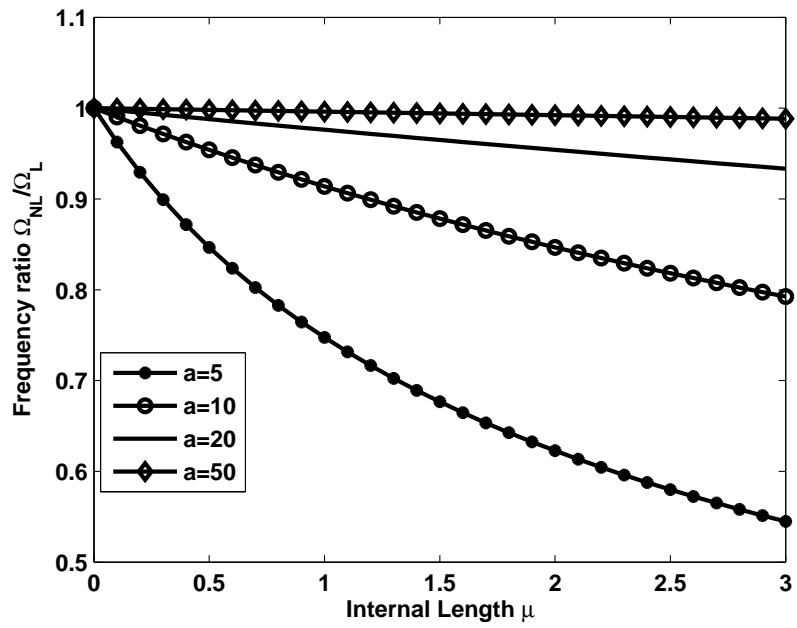}
\caption{Effect of the nonlocal parameter on the frequency ratio $(\Omega_{NL}/\Omega_L$ for a simply supported square plate for various length of the plate with $a/h=$ 10 and gradient index $n=$ 5.}
\label{fig:variouslen}
\end{figure}

The influence of plate dimension on the frequency ratio for various nonlocal parameter is shown in \fref{fig:platedimInfluence} for square plate with plate thickness $a/h=$ 10 and gradient index $n=$ 5. In this case, the nonlocal parameter is assumed to be in the range $\mu=$ (0,3). It can be seen that as the length of the plate increases, the frequency ratio tend to increase and approach the local elasticity solution for considerably larger plate length. From \fref{fig:diffmodes}, it can be seen that the nonlocal parameter has a greater influence on the mode 2 frequency ratio when compared to the mode 1. \fref{fig:variouslen} shows the influence of the nonlocal parameter on the frequency ratio for different plate dimensions with $a/h=$ 10 and gradient index $n=$ 5. It can be seen that as internal length decreases, the frequency ratio increases and approaches the local elasticity solution, irrespective of the plate dimensions. It can be concluded that the effect of the nonlocal parameter is prominent for very thin plates.

% effect of plate length

% effect of internal length on mode1 and mode2

% effect of internal length

\begin{figure}[htpb]
\centering
\includegraphics[scale=0.8]{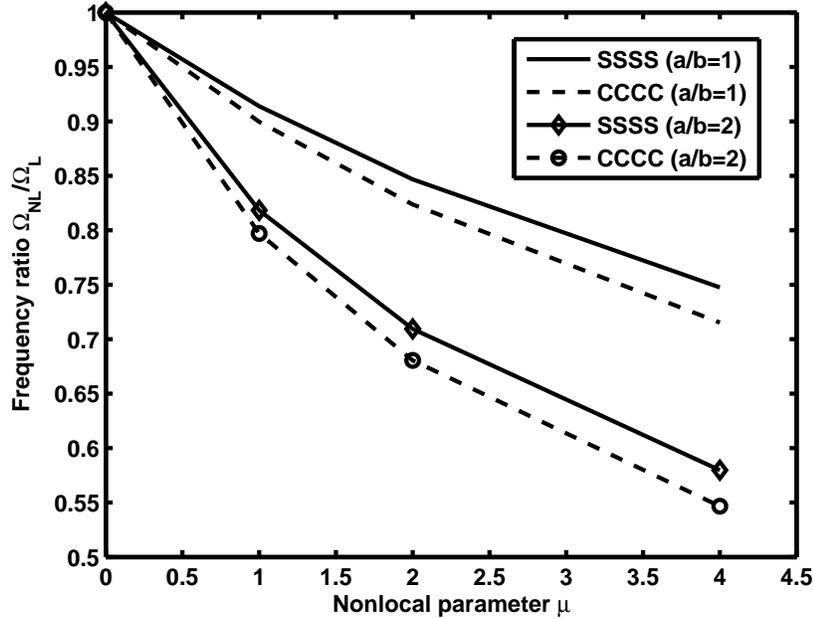}
\caption{Effect of boundary conditions on the frequency ratio $(\Omega_{NL}/\Omega_L$ for a simply supported square and rectangular plate for various internal length with $a/h=$ 10 and gradient index $n=$ 5.}
\label{fig:bceffect}
\end{figure}

% effect of boundary conditions...
The effect of different boundary conditions, viz., all edges simply supported and all edges clamped boundary conditions are graphically shown in \fref{fig:bceffect} and tabulated in Table \ref{table:boundaryeffect} for square and rectangular plate with plate thickness $a/h=$ 10 and gradient index $n=$ 5. It can be seen that both the boundary condition and the nonlocal parameter has an influence on the frequency ratio. The higher the nonlocal parameter, the larger is the influence, irrespective of the boundary condition.

% nondimensionalized frequencies
\begin{table}[htpb]
\renewcommand\arraystretch{1.5}
\caption{Effect of the plate aspect ratio $a/b$, the thickness of the plate $a/h$, the nonlocal parameter $\mu$ and the boundary condition on the natural frequency of a FG plate with $a=$ 10 and gradient index $n=$ 5.} 
\centering
\begin{tabular}{ccclllllll}
\hline
$a/b$ & $a/h$ & $\mu$ & \multicolumn{3}{c}{SSSS} & & \multicolumn{3}{c}{CCCC}\\
\cline{4-10}
& & & mode 1 & mode 2 & mode 3 && mode 1 & mode 2 & mode 3 \\
\hline
\multirow{8}{*}{1} & \multirow{4}{*}{10} & 0 & 0.0441 & 0.1051 & 0.1051 & & 0.0758 & 0.1442 & 0.1455\\
& & 1 & 0.0403 & 0.0860 & 0.0860 & & 0.0682 & 0.1157 & 0.1166 \\
& & 2 & 0.0374 & 0.0745 & 0.0746 & & 0.0624 & 0.0992 & 0.0999 \\
& & 4 & 0.0330 & 0.0609 & 0.0610 & & 0.0542 & 0.0801 & 0.0806 \\
\cline{3-10}
& \multirow{4}{*}{20} & 0 & 0.0113 & 0.0278 & 0.0279 & & 0.0207 & 	0.0410 & 0.0421 \\
& & 1 & 0.0103 & 0.0228 & 0.0228 & & 0.0186 & 0.0326 & 0.0334 \\
& & 2 & 0.0096 & 0.0197 & 0.0198 & & 0.0170 & 0.0279 & 0.0285 \\
& & 4 & 0.0085 & 0.0161 & 0.0162 & & 0.0147 & 0.0225 & 0.0229 \\
\cline{2-10}
\multirow{8}{*}{2} & \multirow{4}{*}{10} & 0 & 0.1055 & 0.1615 & 0.2430 & & 0.1789 & 0.2251 & 0.3022 \\
& & 1 & 0.0863 & 0.1208 & 0.1637 & & 0.1426 & 0.1640 & 0.1960 \\
& & 2 & 0.0748 & 0.1006 & 0.1310 & & 0.1218 & 0.1350 & 0.1557 \\
& & 4 & 0.0612 & 0.0793 & 0.0999 & & 0.0978 & 0.1051 & 0.1181 \\
\cline{3-10}
& \multirow{4}{*}{20} & 0 & 0.0279 & 0.0440 & 0.0701 & & 0.0534 & 0.0682 & 0.0939 \\
& & 1 & 0.0229 & 0.0329 & 0.0464 & & 0.0422 & 0.0491 & 0.0601 \\
& & 2 & 0.0198 & 0.0274 & 0.0371 & & 0.0358 & 0.0402 & 0.0476 \\
& & 4 & 0.0162 & 0.0216 & 0.0283 & & 0.0287 & 0.0313 & 0.0360 \\
\hline
\end{tabular}
\label{table:boundaryeffect}
\end{table}

%--------------- Conclusions
\section{Conclusion}
The linear free flexural vibrations of FG size-dependent nanoplate is numerically studied using NURBS basis functions within a finite element framework. The formulation is based on FSDT and the material is assumed to be graded only in the thickness direction according to the power-law distribution in terms of volume fraction of its constituents. Numerical experiments have been conducted to bring out the effect of the gradient index, the nonlocal parameter, the plate aspect ratio, thickness of the plate and boundary condition on the natural frequency of the FG nanoplate. From the detailed numerical study, it can be concluded that, the natural frequency decreases with increase in the nonlocal parameter and the gradient index. This observation is true, irrespective of thick or thin plate, and square or rectangular plate considered in the present study.

%------------------ bibliography
\bibliographystyle{elsarticle-num}
\bibliography{nle_reference}

\end{document}